\documentclass[12pt]{article}
\usepackage{amsmath}
\usepackage{amssymb}
\usepackage{amscd}
\usepackage{epsfig}
\usepackage{amsthm}

\newcommand{\mob}{\operatorname{Mob}(S^4)}

\newtheorem{dfn}{Definition}[section]
\newtheorem{rem}[dfn]{Remark}
\newtheorem{thm}[dfn]{Theorem}
\newtheorem{defn}[dfn]{Definition}

\newtheorem{prop}[dfn]{Proposition}

\newtheorem{cor}[dfn]{Corollary}
\newtheorem{conj}[dfn]{Conjecture}

\newtheorem{example}[dfn]{Example}

\def\proof{\par\medskip\noindent{\it Proof: }}

\def\C{{\mathbb C}}
\def\R{{\mathbb R}}
\def\H{{\mathbb H}}
\def\Z{{\mathbb Z}}

\def\P{{\mathbb P}}

\def\N{{\mathbb N}}

\def\al{\alpha}
\def\be{\beta}
\def\ga{\gamma}

\def\Ga{\Gamma}

\def\Si{\Sigma}

\def\Om{\Omega}

\def\>{\rangle}
\def\<{\langle}

\def\D{\partial}
\def\3{\ss}
\def\8{\infty}

\def\ol{\overline}

\begin{document}

\title{Conformally flat metrics on 4-manifolds}
\author{Michael Kapovich}
\date{September 28, 2002}

\maketitle

\begin{abstract}
\noindent We prove that for each closed smooth spin 4-manifold $M$ there
 exists a closed smooth 4-manifold $N$ such that $M \# N$ admits a conformally 
flat Riemannian metric.  
\end{abstract}

\tableofcontents

\section{Introduction}

The goal of this note  is to prove 

\begin{thm}
\label{mainthm}
Let $M^4$ be an closed connected smooth spin 4-manifold. 
Then there exists a closed orientable 
4-manifold $N$ such that $M \# N$ admits a conformally 
flat Riemannian metric.  
\end{thm}

Our motivation comes from the following beautiful theorem of C.~Taubes \cite{Taubes}:

\begin{thm}
Let $M$ be a smooth closed oriented 4-manifold. Then there exists a number $k$ 
so that the connected sum of $M$ with $k$ copies of $\ol{\C\P^2}$ 
admits a half-conformally flat structure. 
\end{thm}

Here $\ol{\C\P^2}$ is the complex-projective plane with the 
reversed orientation. Recall that a Riemannian metric $g$ on $M$ is 
{\em anti self-dual} (or half-conformally flat) if the 
self-dual part $W_+$ of the Weyl tensor vanishes. Vanishing of both 
self dual and anti self-dual parts of the Weyl tensor (i.e., vanishing of 
the entire Weyl tensor) is equivalent to local conformal flatness 
of the metric $g$. 

Note that the assumption that $M$ is spin is equivalent to vanishing  
of all Stiefel-Whitney classes, which in turn is equivalent to triviality of the tangent 
bundle of $M'=M\setminus \{p\}$. According to the Hirsch-Smale theory 
(see for instance \cite{Phillips}), 
$M':=M\setminus \{p\}$ is parallelizable iff $M'$ admits an immersion into 
$\R^4$. Thus, by taking $M$ to be simply-connected with nontrivial  
2-nd Stiefel-Whitney class, one sees that $M \# N$ does not admit a flat 
conformal structure for {\em any $N$}: 
otherwise the developing map would immerse  $M'$ into 
$\R^4$. Therefore the vanishing condition is, to some extent, necessary. 
Note also that (unlike in Taubes' theorem) one cannot expect $N$ to be 
simply-connected since the only closed conformally flat simply-connected 
Riemannian manifold is the sphere with the standard conformal structure. 

Sonjong Hwang in his thesis \cite{hwang}, has proven a 3-dimensional version of 
Theorem \ref{mainthm}; moreover, he proves that 
one can use a connected sum of Haken manifolds as the manifold $N$. 
Similar arguments can be used to prove an analogous theorem in the context 
of locally spherical CR structures on 3-manifolds. 

The arguments in both 3-dimensional and 4-dimensional cases, in spirit 
(although, not in the technique), are parallel to Taubes': 
one starts with a singular conformally-flat metric on $M$, where the singularity 
is localized in a ball $B\subset M$. The singular metric is obtained by 
pull-back of the standard metric on the 4-sphere under a branched covering 
$M\to S^4$. Then, by attaching another manifold 
to $M$ along the singular locus (which accounts for taking 
the connected sum $M\# N$), one ``resolves the singularity'' 
and constructs the desired flat conformal structure. 

\begin{conj}
Let $M^n$ be a closed connected smooth n-manifold, whose   
Stiefel-Whitney classes all vanish. Then there exists 
a smooth closed orientable n-manifold $N$ such that $M \# N$ admits a conformally 
flat Riemannian metric.  
\end{conj}

It seems very likely that this conjecture can be proven by methods analogous to the ones 
used in the present paper. The main technical problem is that the reflection groups 
used in the present paper do not exists in the higher dimensions. However, 
it is plausible that one can use instead arithmetic groups containing large 
number of reflections, as it is done in \cite{davis-charney} in a   
different context. 

\medskip 
{\bf Acknowledgments.} During the work on this paper the author was  
 supported by the NSF grants DMS-99-71404 and DMS-02-03045.

\section{Definitions and notation}

We let $\mob$ denote the full group of Moebius transformations of $S^4$, i.e. 
the group generated by inversions in round spheres. Equivalently, $\mob$ is the restriction 
of the full group of isometries $Isom(\H^5)$ to the  $4$-sphere $S^4$ which is the ideal boundary 
of $\H^5$. We will regard $S^4$ as 1-point compactification $\R^4\cup \{\infty\}$ 
of then Euclidean 4-space. 

\begin{defn}
\label{inversion}
Let $Q$ be a unit cube in $\R^4$. We define the {\em PL inversion} $J$ in the boundary of 
$Q$ as follows. Let $h: S^4\to S^4$ be a PL homeomorphism which sends $\Si=\D Q$ onto the round 
sphere $S^3\subset \R^4$ and $h(\infty)=\infty$. Let $j: S^4\to S^4$ be the ordinary inversion 
in $S^3$. Then $J:= h^{-1} \circ j \circ h$. 
\end{defn}

\begin{defn}
A {\em Moebius} or a {\em flat conformal structure} on a smooth $4$-manifold $M$ is an atlas 
$\{ (V_\al, \varphi_\al), \al \in A\}$ which consist of diffeomorphisms 
$\varphi_\al: V_\al \to U_\al \subset S^4$ so that the 
transition mappings $\varphi_\al \circ \varphi_\be^{-1}$ 
are restrictions of Moebius transformations.  
\end{defn}

Equivalently, one can describe Moebius structures on $M$ are conformal classes 
of conformally-Euclidean Riemannian metrics on $M$. Each conformal structure on $M$ gives rise 
to a local conformal 
diffeomorphism, called a {\em developing map}, $d: \tilde{M}\to S^4$, where $\tilde{M}$ 
is the universal cover of $M$. If $M$ is connected, the mapping $d$ is equivariant with respect to a {\em holonomy 
representation} $\rho: \pi_1(M)\to \mob$, where $\pi_1(M)$ acts on $\tilde{M}$ 
as the group of deck-transformations. Given a pair $(d, \rho)$, where $\rho$ is a
representation of $\pi_1(M)$ into $\mob$ and $d$ is a $\rho$-equivariant local diffeomorphism 
from $\tilde{M}$ to $S^4$, one constructs the corresponding Moebius structure on $M$ by taking 
a pull-back of the standard flat conformal structure on $S^4$ to $\tilde{M}$ 
via $d$ and then projecting  the structure to $M$. 

Analogously, one defines a {\em complex-projective structures} on complex 3-manifold $Z$: 
it is a  $\C\P^3$-valued holomorphic atlas on $Z$ so that the transition mappings 
belong to $PGL(3, \C)$. 

\medskip 
The concept of Moebius structure generalizes naturally to the category of orbifolds:  

\noindent A  4-dimensional {\em Moebius orbifold} $O$ is a pair $(X, {\mathcal A} )$, where $X$ is 
a Hausdorff topological space, the {\em underlying space of the orbifold}, 
${\mathcal A}$ is a family of local parameterizations $\psi_\al: U_\al\to U_\al/\Ga_\al=V_\al$, 
where $\{V_\al, \al\in A\}$ is an open covering of $X$, 
$U_\al$ are open subsets in $S^4$, $\Ga_\al$ are finite groups of Moebius automorphisms  
of $U_\al$ and the mappings $\psi_\al$ satisfy the usual compatibility conditions:

If $V_\al\to V_\be$ is the inclusion map then we have a Moebius embedding   
$U_\al\to U_\be$ which is equivariant with respect to a monomorphism $\Ga_\al \to \Ga_\be$, so that 
the diagram
$$
\begin{array}{ccc}
U_\al & \to & U_\be\\
\downarrow & ~ & \downarrow\\
V_\al & \to & V_\be
\end{array}
$$ 
is commutative. The groups $\Ga_\al$ are the {\em local fundamental groups} of the orbifold $O$. 

\medskip
For a Moebius orbifold one defines a developing mapping $d: \tilde{O}\to S^4$ 
(which is a local homeomorphism 
from the universal cover $\tilde{O}$ of $O$) and, if $O$ is connected, 
a holonomy homeomorphism $\rho: \pi_1(O)\to \mob$, 
which satisfy the same equivariance condition as in the manifold case. Again, given a 
pair $(d, \rho)$, where $d$ is a $\rho$-equivariant homeomorphism, one defines 
the corresponding Moebius structure via pull-back. 

\begin{example}
Let $G\subset \mob$ be a subgroup acting properly discontinuously on an open subset $\Om\subset S^4$. 
Then quotient space $\Om/G$ has a natural Moebius orbifold structure. The local charts $\phi_\al$ 
appear in this case as restrictions of the projection $p: \Om\to \Om/G$ to open subsets 
with finite stabilizers. 
\end{example}

In particular, suppose that $G$ is a finite subgroup of $\mob$ generated by reflections, the quotient 
$Q:= \Om/G$ can be identified with the intersection of a fundamental domain of $G$ 
with $\Om$. The Moebius structures on 4-dimensional manifolds and orbifolds constructed 
in this paper definitely do not arise this way. 
A more interesting example is obtained by taking a manifold $M$ and a local 
homeomorphism $h: M\to S^4$, so that $Q\subset h(M)$. Then we can pull-back the Moebius 
orbifold structure on $Q$ to an appropriate subset $X$ 
of $M$, to get a 4-dimensional Moebius orbifold. As another example of a pull-back construction, let 
 $O$ be a Moebius orbifold and $M\to O$ be an orbifold cover such that $M$ is a manifold. Then 
 one can pull-back the Moebius orbifold structure from $O$ to a usual Moebius structure on $M$.

\section{Reflection groups in $S^4$ with prescribed combinatorics of the fundamental domains}
\label{reflect}

\begin{figure}[tbh]
\centerline{\epsfxsize=3.5in \epsfbox{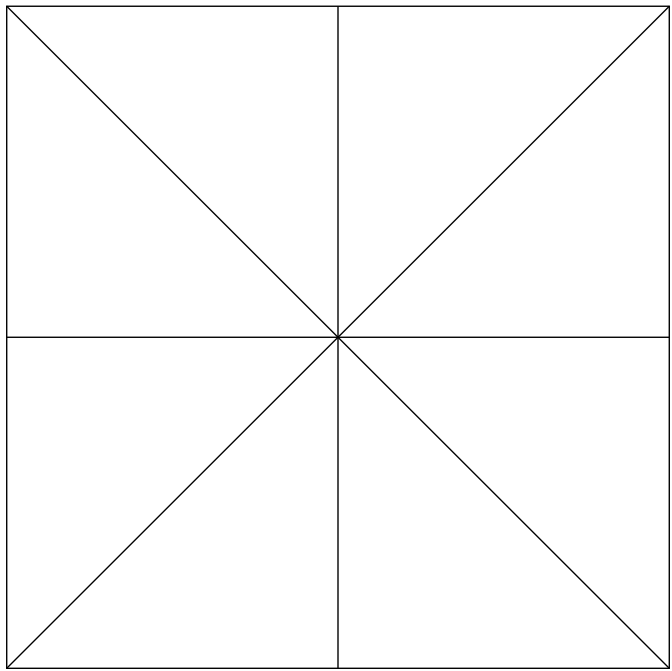}}
\caption{\sl Barycentric subdivision of a square.}
\label{F2}
\end{figure}

Consider the standard cubulation ${\mathcal Q}$ of $\R^4$ by the unit Euclidean cubes and let $X$ denote the 
2-skeleton of this cubulation. Given a collection of round balls $\{B_i, i\in I\}$ in $\R^4$, 
with the nerve ${\mathcal N}$, 
we define the canonical simplicial mapping $f: {\mathcal N}\to \R^4$ by sending each vertex of 
${\mathcal N}$ to the center of the corresponding ball and extending $f$ linearly to the simplices 
of ${\mathcal N}$. For a subcomplex $K\subset {\mathcal Q}$ define its {\em barycentric subdivision} 
$\beta(K)$ to be the following simplicial complex. Subdivide each edge of $K$ by its midpoint. 
Then inductively subdivide each $k$-cube $Q$ in $K$ by coning off the barycentric subdivision of $\partial Q$  
from the center of $Q$.

\begin{figure}[tbh]
\centerline{\epsfxsize=3.5in \epsfbox{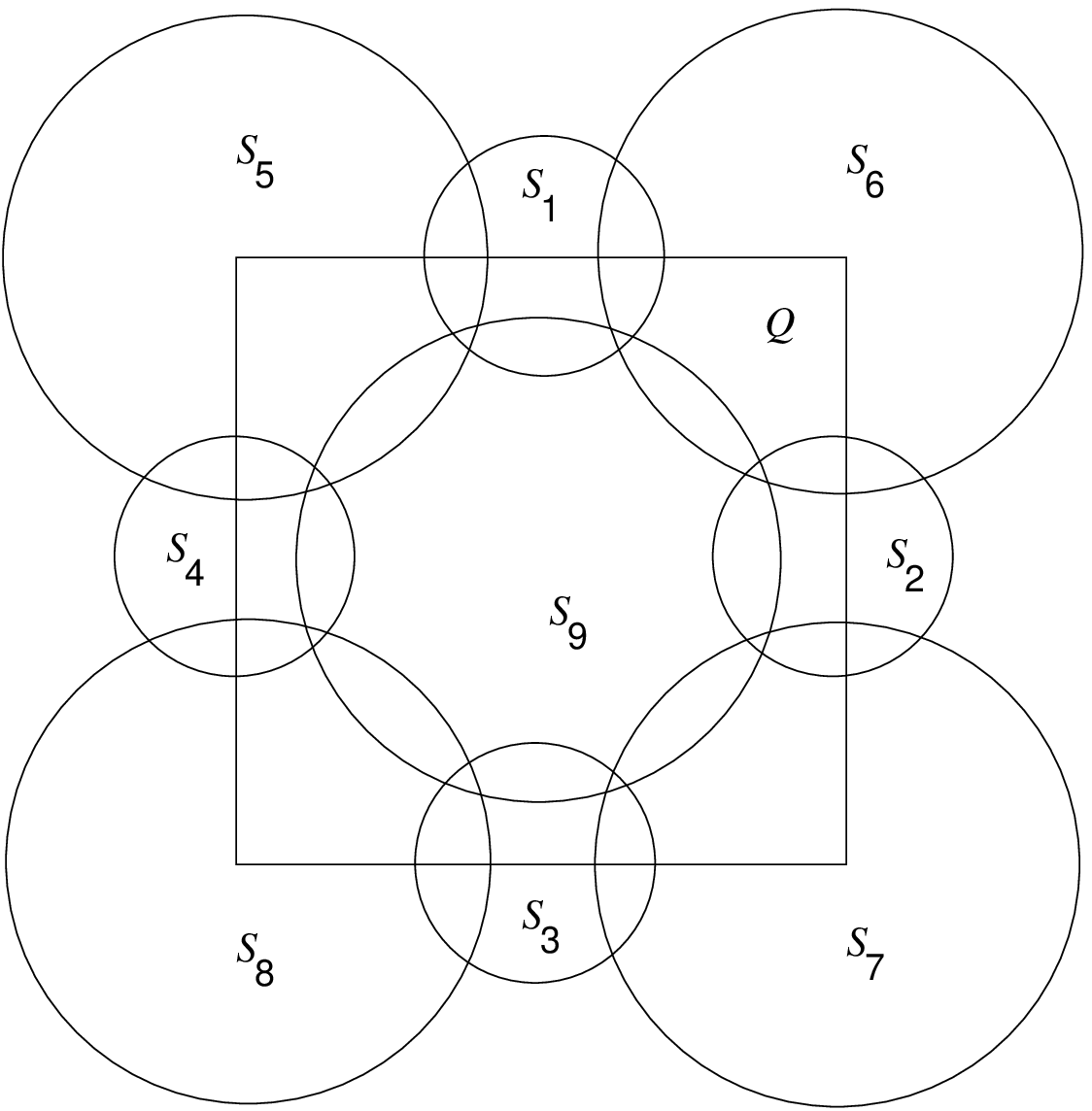}}
\caption{}
\label{F1}
\end{figure}

\begin{prop}
Suppose that $K\subset X$ is a 2-dimensional compact subcomplex such that each vertex belongs to 
a 2-cell. Then there exists a collection of open round 4-balls $B_i$, $i=1,...,k$, centered 
at the vertices of $\beta(K)$, so that: 

(1) The Moebius inversions $R_i$ in the round spheres $S_i=\D B_i$ generate a 
discrete reflection group $G\subset \mob$.

(2) The complement $S^4 \setminus \cup_{i=1}^k B_i$ is a fundamental domain $\Phi$ of $G$. 

(3)  The canonical mapping from the nerve of $\{B_i, i=1,...,k\}$ to $\R^4$ is a simplicial isomorphism  
onto $\beta(K)$.  
\end{prop}
\proof We begin by constructing the family of spheres $S_i, i\in \N$ centered at certain points of $X$. 
For each square $Q$ in $X$ we pick 9 points $x_1,....,x_9$: $x_5,...,x_8$ are the vertices of $Q$,    
$x_1,...,x_4$ are midpoints of the edges of $Q$ and $x_9$ is the center of $Q$. 
We then take spheres $S(x_i, r)$ of radius $r=\frac{1}{\sqrt{6}}$ centered at the points $x_4,...,x_9$ 
and the spheres $S(x_i, \rho)$ of radius $\rho= \frac{1}{\sqrt{12}}$ centered at the points $x_1,...,x_4$. 
See Figure \ref{F1}. The reader will verify that:

1. The spheres $S(x_i, r)$ and $S(x_j, r)$ are disjoint provided that 
$i\ne j\in \{1,...,4\}$ and $i\ne j\in \{5,...,8\}$. 

2. The spheres $S(x_i, r)$ and $S(x_9, r)$ intersect at the right angle, $i=1,...,4$;   
the spheres $S(x_i, r)$ and $S(x_9, r)$ intersect at the (exterior) angle $\frac{\pi}{3}$, $i=5,...,8$. 

3. The spheres $S(x_i, r)$ and $S(x_j, r)$ intersect at the (exterior) angle $\frac{\pi}{3}$, 
provided that $1\le i\le 4$, $1\le j\le 8$ and $j=5+i, 6+i$, and are disjoint otherwise. 

\medskip
Suppose now that $Q^4$ is a unit 4-cube, apply the above construction to each 2-face of $Q^4$. 
The reader will verify that the properties 1--3 of the spheres $S_i$,  
ensure that the covering $\{B_i\}$ of $(Q^4)^{(2)}$ has the nerve ${\mathcal N}_{Q^4}$ such that the canonical mapping 
${\mathcal N}_{Q^4}\to \beta((Q^4)^{(2)})$ is a simplicial isomorphism. 

Now we are ready to construct the covering $\{B_i: i=1,...,k\}$ of the 2-complex $K$. For each 2-face $Q$ of $K$ introduce the 
family of nine round spheres $S_i$ constructed above, consider the inversions 
$R_i$ is these spheres; the spheres $S_i$ bound balls $\{B_i: i=1,...,k\}$. 
The fact that for each 4-cube $Q^4$ the mapping 
${\mathcal N}_{Q^4}\to \beta((Q^4)^{(2)})$ is a simplicial isomorphism, implies that the mapping 
from the nerve of the covering  $\{B_i: i=1,...,k\}$ to $K$ is a simplicial isomorphism as well. 
Thus the exterior angles of intersections between the spheres equal $\frac{\pi}{2}$ and  
$\frac{\pi}{3}$, thus we can apply Poincare's fundamental polyhedron theorem \cite{Maskit} to ensure that the 
intersection of the complements to the balls $B_i$ is a fundamental domain for 
the Moebius group $G$ generated by the above reflections.  
 \qed

\begin{rem}
Instead of collections of round balls based on a cubulation of $\R^4$ one can use a periodic triangulation 
of $\R^4$, however in this case the construction of a collection of balls covering the 2-skeleton of 
a 4-simplex is slightly more complicated. 
\end{rem}

\section{Proof of Theorem \ref{mainthm}}

Recall that vanishing of all Stiefel-Whitney classes of the manifold $M$ implies that the manifold 
$M'= M\setminus \{p\}$ is parallelizable; hence, by \cite{Phillips}, there exists an immersion 
$f: M'\to \R^4$. Let $B$ denote a small open round ball centered at $p$ and let 
$M''$ denote the complement $M\setminus B$. We retain the notation $f$ for the restriction $f|M''$.  
We next convert to the piecewise-linear setting, since in dimension 4 the categories of 
PL and smooth structures are the same it does not limit our discussion: 
cubulate the manifold $M$ so that $\D B$ is a subcomplex of the cubulation 
and that the restriction of $f$ to each 4-cube $Q'\subset M''$ 
is a diffeomorphism onto a cube in the standard (unit cube) cubulation of $\R^4$. 
Without loss of generality we may assume that the mapping $f$ preserves the orientation. 

We now borrow the standard arguments from the proof of Alexander's theorem which states that each 
closed $n$-dimensional PL manifold is a branched cover over the $n$-sphere, see e.g. \cite{Feighn}.  
Extend the map $f$ to a map $F$ on the ball $B$ so that 
the restriction of $F$ to each 4-cube is a diffeomorphism onto a cube in the standard cubulation of $\R^4$. 
Now, for each cube $Q'\subset M$ such that $F|Q'$ is orientation-reversing we replace 
$F|Q'$ with the composition $J\circ F|Q'$, where $J$ is the PL {\em inversion} in the 
boundary of the unit cube $F(\D Q')$ (see Definition \ref{inversion}). 
The resulting mapping $h: M\to S^4$ has the property that it is a local PL homeomorphism 
away from a 2-dimensional subcomplex $L\subset B$. (Note that $L$ has dimension 2 
near every point: each vertex in $L$ belongs to a 2-cube.) Thus the mapping $h$ is a branched 
covering over $S^4$ with the singular locus $L$ contained in the ball $B$, the branch-locus 
of $h$ is the compact subcomplex $K=h(L)\subset \R^4$.   The branched covering $h$ has the 
property that for each point $x\in K$ there exists a neighborhood $U(x)\subset \R^4$ 
such that $h^{-1}(U(x))$ is a disjoint union of balls $V(y), y\in h^{-1}(x)\in B$, 
(whose interiors contain $y$), so that for each $y\in h^{-1}(x)$, the restriction  
$h|V(y)$ is a branched covering onto $U(x)$. Moreover, each branched 
covering $h|V(y)$ is obtained by coning off a branched covering from the 3-sphere $\D V(y)$ 
to the 3-sphere $U(x)$. 

Let $T$ denote a regular neighborhood of $K$ in $\R^4$, so that $U(x)\subset T$ for each $x\in K$. 
Next, subdivide the cubulation of $\R^4$ and scale the subdivision up to the standard unit cubulation, so that 
the discrete group $G$ and the collection of balls $\{ B_j, j=1,...,k\}$ associated with the subcomplex $K$ 
in section \ref{reflect} have the properties: 

1. $T\subset \cup_{i=1}^k B_k$. 

2. Each ball $B_j$, $j=1,...,k$, (centered at $x_j\in K$) is contained in the neighborhood  
$U(x_j)$. 

\medskip 
We now use the branched covering $h$ to introduce a Moebius orbifold structure $O$ on the 
complement $X_O$ to an open tubular neighborhood ${\mathcal N}^0(L)$ of $L$ in $M$ as follows:

For each ball $B_j\subset U(x_j)$ centered at $x_j\in K$ and for each $y_j\in h^{-1}(x_j)\cap L$, such 
that the restriction $h|V(y_j)$ is {\em not a homeomorphism onto its image}, we let $\tilde{B}(y_j)$ 
denote the inverse image $h^{-1}(B_j)\cap V(y_j)$. It follows that each $\tilde{B}(y_j)$  
is a polyhedral 4-ball in $M$ and the union of these balls is a tubular neighborhood ${\mathcal N}(L)$ of $L$. 
The boundary of ${\mathcal N}(L)$ has a natural partition into subcomplexes: ``vertices'', ``edges'', 
``2-faces'' and ``3-faces'':

\begin{itemize}
\item The ``vertices'' are the points of triple intersections of the 3-spheres 
$\partial \tilde{B}(y_j)$, $\partial\tilde{B}(y_i)$, $\partial\tilde{B}(y_l)$.  

\item The ``2-faces'' are the connected components of the double intersections of 
the 3-spheres $\partial \tilde{B}(y_j)$, $\partial\tilde{B}(y_i)$.  

\item The ``3-faces'' are the connected components of the complements 
$$
\partial \tilde{B}(y_j) \setminus \cup_{i\ne j} \tilde{B}(y_i).$$ 
\end{itemize}

We declare each ``3-face'' a boundary  reflector of the orbifold $O$. The dihedral angles 
between the balls $B_j$ define the dihedral angles between the boundary reflectors in $O$. 
Since the restriction $h| M \setminus L$ is a local homeomorphism, this construction defines 
a Moebius orbifold $O$. The mapping $h|X_O$ is the projection of the developing mapping 
$\tilde{h}: \tilde{O}\to S^4$ of this Moebius orbifold. Let $O'$ denote the  
orbifold with boundary $B\cap O$; let $Q$ be the closed orbifold obtained by 
attaching 4-disk $D^4$ along the boundary sphere $S^3=\partial B$.  

We now convert back to the smooth category. 
It is clear from the construction that the orbifold $O$ is obtained by (smooth) gluing of the 
manifold with boundary $M\setminus B$ and the orbifold with boundary $O'$. Hence 
$O$ is diffeomorphic to the connected sum of the manifold $M$ with the orbifold $Q$. 
We also note that all local fundamental groups of $O$ embed naturally into $\pi_1(O)$. 

It remains to construct a finite manifold covering $\hat{M}$ over the orbifold $O$, so that 
$M\setminus B$ lifts homeomorphically to $\hat{M}$; the construction is analogous to the 
one used by M.~Davis in \cite{davis}. The fundamental group $\pi_1(O)$ 
is the free product $\pi_1(M)*\pi_1(Q)$. We have holonomy homomorphism 
$$
\phi: \pi_1(O)\to G,
$$
the subgroup $\pi_1(M)$ is contained in the kernel of this homomorphism; by construction, 
the kernel of $\phi$ contains no elements of finite order. The Coxeter group $G$ is virtually torsion-free, 
let $\theta: G\to A$ be a homomorphism onto a finite group $A$, so that $Ker(\theta)$ is torsion-free and 
orientation-preserving. 
Then the kernel of the homomorphism $\psi=\theta\circ \phi: \pi_1(O)\to A$ is a torsion-free finite 
index subgroup of $\pi_1(O)$, which contains $\pi_1(M)$. Let  $\hat{M}\to O$ denote the 
finite orbifold cover corresponding to the  subgroup $Ker(\psi)$. Then $\hat{M}$ is a smooth oriented 
conformally flat manifold, the submanifold $M\setminus B$ lifts diffeomorphically into $M\setminus B\subset \hat{M}$. 
Thus the connected sum decomposition $O= M\# Q$ also lifts to 
$\hat{M}$, so that the latter manifold is diffeomorphic to the connected sum of $M$ and a 4-manifold $N$. 
\qed 

\bigskip 
We observe that the proof of Theorem \ref{mainthm} can be modified to prove the following: 

\begin{thm}
Suppose that $M$ is a closed smooth 4-manifold with vanishing second Stiefel-Whitney class. Then 
there exists a closed smooth 4-manifold $N$ so that $\hat{M}=M\# N$ admits a 
conformally-Euclidean Riemannian metric. 
\end{thm}
\proof The difference with Theorem \ref{mainthm} is that $M$ can be nonorientable. Let 
$\tilde{M}\to M$ be the orientable double cover with the deck-transformation group $D\cong \Z/2$. 
Then all Stiefel-Whitney classes of $\tilde{M}$ are trivial. 
As before, let $p\in M$, $\{p_1, p_2\}$ be the preimage of $\{p\}$ in $M$. Consider 
a Euclidean reflection $\tau$ in $\R^4$ and an epimorphism $\theta: D\to \<\tau\>$. Then,  
arguing as in the proof of Phillips' theorem \cite{Phillips}, one gets a $\theta$-equivariant immersion 
$\tilde{f}: \tilde{M}\setminus \{p_1, p_2\} \to \R^4$. This yields a $D$-invariant flat conformal structure 
on $\tilde{M}\setminus \{p_1, p_2\}$ via pull-back of the flat conformal structure from $\R^4$. 
Let $B_1\sqcup B_2$ be a $D$-invariant disjoint union of open balls around the points $p_1, p_2$.   
Then the rest of the proof of Theorem \ref{mainthm} goes through: replace the ball $B_1$ 
with a manifold with boundary $N_1$ so that the flat conformal structure on  
$\tilde{M}\setminus (B_1\cup B_2)$ extends over $N_1$. Then glue a copy of $N_1$ along the boundary 
of $B_2$ in $D$-invariant fashion. Note that the quotient of 
the manifold $P:=(\tilde M\setminus (B_1 \cup B_2))\cup (N_1 \cup N_2)$ by the group $D$ is diffeomorphic 
to a closed manifold $M\# N$, where $N$ is obtained from $N_1$ by attaching the 4-ball along the boundary.  
Finally, project the $D$-invariant Moebius structure on $P$ to a Moebius structure on the manifold $M\# N$. 
\qed

As  a corollary of Theorem \ref{mainthm} we get:

\begin{cor}
Let $\Ga$ be a finitely-presented group. Then there exists a 3-dimensional complex manifold $Z$ 
which admits a complex-projective structure, so that the fundamental group of $Z$ splits as 
$\Ga * \Ga'$. 
\end{cor}
\proof Our argument is similar to the one used to construct (via Taubes' theorem) 
3-dimensional complex manifolds with the prescribed finitely-presented fundamental group, see \cite{Taubes}.  
We first construct a smooth closed oriented 4-dimensional spin 
manifold $M$ with the  fundamental group $\Ga$. This can be done for instance as follows. 
Let $\< x_1,...,x_n | R_1,..., R_\ell\>$ be a presentation of $\Ga$. 
Consider a 4-manifold $X$ which is the connected sum of $n$ copies of $S^3\times S^1$. 
This manifold is clearly spin. Pick a collection of disjoint embedded smooth 
loops $\ga_1,...,\ga_\ell$ in $X$, which represent the conjugacy classes of the words 
$R_1,...,R_\ell$ in the free group $\pi_1(X)$. Consider the pair $(S^4, \ga)$, where $\ga$ 
is an embedded smooth loop in $S^4$. For each $i$ pick a diffeomorphism $f_i$ between 
a tubular neighborhood $T(\ga)$ of $\ga$ in $S^4$ and a tubular neighborhood $T(\ga_i)$ of $\ga_i$ in $X$. 
We can choose $f_i$ so that it matches the spin structures of $T(\ga)$ and 
$T(\ga_i)$. Now, attach $n$ copies of $S^4\setminus T(\ga)$ to $X\setminus \cup_i T(\ga_i)$ 
via the diffeomorphisms $f_i$. 
The result is a smooth spin 4-manifold $M$ with the fundamental group $\Ga$. 

Next, by Theorem \ref{mainthm} there exists 
a smooth 4-manifold $N$ (with the fundamental group $\Ga'$) 
such that $\hat{M}=M\# N$ admits a conformally-Euclidean Riemannian metric. 
Applying the twistor construction to the manifold $\hat{M}$ we get a complex 3-manifold $Z$ which 
is an $S^2$-bundle over $\hat{M}$ and the flat conformal structure on $\hat{M}$ lifts 
to a complex-projective structure on $Z$, see for instance \cite{twist}. 
Clearly, $\pi_1(Z)\cong \pi_1(\hat{M})= \Ga * \Ga'$. \qed

\label{bib}

\noindent Department of Mathematics, University of Utah,\\
Salt Lake City, UT 84112-0090,\\
kapovich$@$math.utah.edu

\end{document}